\documentclass[12pt]{article}
\usepackage{amsmath,amssymb}
\newcommand{\bs}{\bigskip}

\newcommand{\mb}{\mbox}

\newcommand{\stt}{\subseteq}
\newcommand{{\z}}{\mathbb Z}
\renewcommand{\P}{\mathbb P}
\newcommand{\irr}{irreducible }
\newcommand{\CM}{Cohen-Macaulay }
\newcommand{\Hs}{Hilbert scheme }

\newcommand{\alg}{algebraic }
\newcommand{\emb}{embed}
\renewcommand{\c}{curve}

\begin{document}
\setlength{\baselineskip}{16pt}

\begin{center}{\large\bf Questions of Connectedness of}\\
{\large\bf the Hilbert
Scheme of Curves in ${\mathbb P}^3$}

\bs\bs{\sc Robin Hartshorne}\bs\\
Department of Mathematics\\University of California\\
Berkeley, California 94720--3840

\vspace{.5 truein}
{\it Dedicated to S.~Abhyankar on the occasion of his
70th birthday.}\end{center}

\bs\bs
\begin{quote}
We review the present state of the problem, for each degree $d$
and genus $g$, is the Hilbert scheme of locally Cohen--Macaulay
curves in ${\mathbb P}^3$ connected?
\end{quote}

\section{Introduction}

In studying \alg \c s in projective spaces, our forefathers in the 19th
century noted that \c s naturally move in \alg families. In the
projective plane, this is a simple matter. A \c \ of degree $d$ is
defined by a single homogeneous polynomial in the homogeneous
coordinates $x_0,x_1,x_2$. The coefficients of this polynomial give a
point in another projective space, and in this way \c s of degree $d$ in
the plane are parametrized by the points of a ${\P}^N$ with
$N=\frac 12 d(d+3)$. For an open set of ${\P}^N$, the corresponding \c \
is \irr and nonsingular. The remaining points of
${\P}^N$ correspond to \c s that are singular, or reducible, or have
multiple components. In particular, the nonsingular \c s of degree $d$
in ${\P}^2$ form a single \irr family.

In ${\P}^3$, the situation is more complicated. For a given degree $d$,
there may be \c s with several different values of the genus $g$.
Even for fixed $d,g$, the family of \c s with given $d,g$ may not be
\irr\!. An early example, noted by Halphen and Weyr in 1874 is the case
$d=9$ and $g=10$. One type consists of \c s $C_1$ of bidegree (3,6)
on a nonsingular  quadric surface $Q$. The other type consists of \c s
$C_2$ that are the complete intersection of two cubic surfaces $F_3$ and
$F'_3$. These \c s form two \irr components of the \Hs $H^{sm}_{d,g}$
of smooth \c s of degree $d$ and genus $g$ in ${\P}^3$.
Furthermore, it is not hard to see that every \c \ in $H^{sm}_{9,10}$
belongs to one  of these two types, and that there is no flat family of
\c s whose general member belongs to one type and whose special member
belongs to the other type \cite[Ch.IV, 6.5.4]{H2}. Thus the \Hs of
smooth \c s of given degree and genus in ${\P}^3$ need not be connected.

On the other hand, if one phrases the question more generally, by
letting a ``\c\!" mean an arbitrary closed subscheme of dimension 1 in
${\P}^3$, then the \Hs for each degree $d$ and arithmetic genus $g$ is
connected. In fact, in my thesis \cite{H1}, I showed that the \Hs of
closed subschemes of ${\P}^n$ with Hilbert polynomial $P$ is connected
(provided it is nonempty) for any $n$ and any $P$. In the proof,
non-reduced schemes play an essential role. Here is the main idea of the
proof for the case of \c s in ${\P}^3$. Suppose, for example, that we
start with a nonsingular \c \ $C$. Its general projection to ${\P}^2$
will be a plane \c \ $C_0$ with nodes. Using the projection we can
construct a flat family whose general member is $C$ and whose special
member $C_1$ is a \c \ with support $C_0$, having \emb ded points at the
nodes (see \cite[III, 9.8.4]{H2}) for an example showing how these \emb
ded points arise). Then we can make another flat family, pulling the
\emb ded points off $C_1$, to get $C_0$ union a number of points in
${\P}^2$.  Finally, we move $C_0$ in a flat family of plane \c s to a
union of lines in ${\P}^2$ meeting at a single point. If the original \c
\ $C$ had degree $d$ and genus $g$, we obtain in this way a ``fan" of
$d$ lines in the plane together with $k=\frac 12
(d\!-\!1)(d\!-\!2)\!-\!g$ isolated points in the plane. Any other \c \
$C'$ with the same $d,g$ can be connected by a sequence of flat
specializations and generalizations to the same configuration, so $C$
and $C'$ are connected within the \Hs $H_{d,g}$, of all closed
subschemes of ${\P}^3$ of dimension 1, degree $d$, and arithmetic genus
$g$. For nonreduced \c s a slightly more complicated, but similar method
applies.

Thus we have a connectedness theorem for the \Hs of curves in ${\P}^3$,
but it is unsatisfactory in that, even if we want to connect one smooth
\c \ to another, we must pass by way of schemes with \emb ded points and
isolated points, which one can argue should not really count as ``\c s".

With the development of liaison theory in recent years
\cite{MP2},\cite{H5}, an intermediate class of \c s has received
much attention, the locally \CM curves. We say a curve is locally \CM if
it is a scheme of equidimension 1, and all its local rings are \CM
rings. Equivalently it is a 1-dimensional scheme with no \emb ded
points or isolated points.
It is clear that this class of \c s is the natural class in which to do
liaison: even if one is primarily interested in nonsingular \c s, the
minimal \c s in a biliaison class may be reducible and non-reduced. So
we pose the question: Is the \Hs $H^{CM}_{d,g}$ of locally \CM \c  s of
degree $d$ and arithmetic genus $g$ in ${\P}^3$ connected?
The answer is unknown at present, so we devote this paper to a survey of
the current state of this question.

\section{Known results}

\subsection{When is $H_{d,g}$ nonempty?}  Before
discussing whether a given \Hs is connected, one should at least know
when it is nonempty.

For smooth \c s in ${\P}^3$, the result was stated by Halphen \cite{Ha}
with an incorrect proof, and proved one hundred years later by Gruson and
Peskine \cite{GP},\cite{H3}. There exists an \irr smooth curve $C$ of
degree $d$ and  genus $g$ in ${\P}^3$ if and only if either
\begin{quote}
a) \ $d\geq 1$ and $g=\frac 12(d-1)(d-2)$ (these are the plane \c s),
or\\ b) \ there exist $a,b>0$ with $d=a+b$ and $g=(a-1)(b-1)$\\
\mbox{}\quad \ (these are \c s
on quadric surfaces), or\\
c)  \ $d\geq 3$ and $0\leq g\leq \frac 16 d(d-3)+1$.
\end{quote}
The hardest part of the proof is the existence of \c s for all $(d,g)$
in the range c), which they construct on suitable cubic and quartic
surfaces in ${\P}^3$.

If one considers all one-dimensional closed subschemes of ${\P}^3$,
the answer was known to Macaulay \cite{Mac}, and rediscovered in
\cite{H1}. Then $H_{d,g}$ is nonempty for all $d\geq 1$ and all
arithmetic genus $g\leq \frac 12(d-1)(d-2)$. The existence is simple.
Just take a plane \c \
  of degree $d$ and add lots of isolated points. Note that the arithmetic
genus $g$ can become arbitrarily negative.

For locally \CM \c s, the answer is slightly more complicated,
but not too difficult \cite{H4}. A locally \CM \c \ with given $d,g$
exists if and only if either
\begin{quote}
a) $d\geq 1$, \ $g=\frac 12(d-1)(d-2)$ (a plane \c s), or\\
b) $d\geq 2$, \ $g\leq\frac 12(d-2)(d-3)$
\end{quote}
For $d=2$ one can exhibit a multiplicity two structure on a line with
any given arithmetic genus $g\leq 0$. For example, the scheme in ${\P}^3$
defined by the homogeneous ideal $(x^2,xy,y^2,xz^r-yw^r)$, for any
$r\geq 0$, has
$g=-r$.  Then one can construct \c s for all $(d,g)$ in case b) above by
taking a plane \c \ of degree $d-1$ containing a line, and putting a
suitable multiplicity two structure on the line.

\subsection{When is $H_{d,g}$ \irr\!?}
 From now on, we will consider only locally \CM \c s, and denote
$H^{CM}_{d,g}$ by $H_{d,g}$.

There are some values of $(d,g)$ for which $H_{d,g}$ is \irr\!,
and hence trivially connected \cite{MP3}. These are
\begin{quote}
a) \ $d\geq 1$, \ $g=\frac 12(d-1)(d-2)$, the plane \c s\\
b) \ $d=2$, $g\leq 0$. For $g=0$ we have a plane curve; for $g=-1$,
two disjoint lines or a double line on a quadric; and for $g\leq -2$,
double structures on a line.\\
c) \ Some special values of $g$ for higher degree, namely
$(d,g)=(3,0),(3,-1),(4,1)$, and
\[
d\geq 5 \ , \qquad
\textstyle{\frac 12}(d-3)(d-4)+1 < g\leq \textstyle{\frac 12}
(d-2)(d-3) \ .
\]
\end{quote}
For all other $(d,g)$, namely
$d=3$, $g\leq -2$; \ $d=4$, $g\leq 0$; and $d\geq 5$,
$g\leq \frac 12 (d-3)(d-4)+1$, $H_{d,g}$ has two or more \irr components.

\subsection{Extremal \c s}
For any curve $C\stt{\P}^3$, an important invariant is the {\it Rao
module} $M(C)=\bigoplus_{n\in{\z}} H^1({\mathcal I}_C(n))$. The
dimensions of the graded components of the Rao module are the {\it Rao
function} $\rho_C(n)=\dim H^1({\mathcal I}_C(n))$. A \c \ is
arithmetically \CM (ACM) if and only if its Rao module is 0.

For non-plane \c s, there are explicit bounds on the Rao function in
terms of $d$ and $g$ \cite{MP4}. In particular, for all $n$ we have
\[
\rho(n)\leq  \textstyle{\frac 12}
(d-2)(d-3)-g \ .
\]
Thus, if $g= \frac 12 (d-2)(d-3)$, the Rao function is 0, so the \c \ is
necessarily ACM, and one knows in this case that the \Hs is \irr
\cite{E}.

If $g< \frac 12
(d-2)(d-3)$, then one has the more precise result that $\rho(n)$ is
bounded by a function that is equal to $\frac 12
(d-2)(d-3)-g$ for $0\leq n\leq d-2$, and decreases with slope 1
(resp.~--1) to zero on both ends of this range.

In their paper \cite{MP3}, Martin-Deschamps and Perrin define an {\it
extremal} \c \ to be a non-ACM \c \ whose Rao function is equal to this
bound for all $n$. For any $g\!<\! \frac 12
(d\!-\!2)(d\!-\!3)$, they show the existence of extremal \c s, and show
that they form an \irr component of the \Hs\!.

For \c s that are not extremal, Nollet has established a stronger bound
on the Rao function \cite{N1}. If $d\geq 5$ and the \c \ is neither ACM
nor extremal, then $\rho(n)\leq\linebreak\frac 12
(d\!-\!3)(d\!-\!4)+\!1-g$. In particular this implies  $g\leq\frac 12
(d-3)(d-4)+1$. Thus any \c \ of
$d\geq 5$ and $g>\frac 12 (d-3)(d-4)+1$  must be extremal, and we
conclude that the \Hs is \irr in that range. Curves satisfying Nollet's
stronger bounds are called {\it subextremal.}

\subsection{When is $H_{d,g}$ connected?}
If the \Hs has two or more \irr components, which happens for
$d=3$, $g\leq-2$; \ $d=4$, $g\leq 0$; and
$d\geq 5$, $g>\frac 12 (d-3)(d-4)+1$ we can ask if it is connected.
Here are some cases in which it is known to be connected.
\begin{enumerate}
\item \ If $d=3$, $g\leq -2$, then $H_{d,g}$ is connected
\cite{N2}, and has approximately $\frac 13|g|$ \irr components.
\item \ If $d=4$, $g\leq 0$, then $H_{d,g}$ is connected \cite{NS}.
\item \ If $d\geq 5$ and $g>\frac 12 (d-3)(d-4)+1$, then $H_{d,g}$ has
two \irr components, consisting of the extremal \c s in one and the ACM
\c s in the other, and  is connected \cite{N3}.
\item \ If $d\geq 4$ and $g=\frac 12 (d-3)(d-4)$, then $H_{d,g}$ has
2, 3, or 4 \irr components,  and  is connected \cite{A1}.
\item \ If $d=5$, $g=0$, $H_{d,g}$ has four \irr components
  and  is connected \cite{H9}.
\item If $d\geq 6$ and $g=\frac 12 (d-3)(d-4)-1$, then $H_{d,g}$ has
four or five irreducible components, and is connected \cite{Sa}.
\end{enumerate}

These cases, together with the cases of $H_{d,g}$ \irr listed above, are
the only cases in which it is known that $H_{d,g}$ is connected at
present. The problem falls into two halves. The first is to list the
\irr components of $H_{d,g}$, and the second is to show the existence of
flat families of \c s connecting the different components. It is the
first of these that is blocking further progress at the moment, because
it requires a classification of all \c s of the given $d,g$. The most
difficult part is to understand all the possible nonreduced structures
on a \c \ of lesser degree. Thus already the case of multiplicity four
structures on a line is extremely complicated.

To avoid the first problem, we formulate the question differently.

\subsection{Curves connected to extremal \c s}
For a given $d,g$, there is always one \irr component of $H_{d,g}$
consisting of the extremal \c s. So we ask, which classes of \c s can be
connected by flat families in $H_{d,g}$
to an extremal \c\!? If every \c \ with the given $(d,g)$ is connected
to an extremal \c\!, then $H_{d,g}$ is connected.

The advantage of this question is that we do not have to classify all \c
s of type $(d,g)$. Here are some cases that are known, namely \c s that
can be connected within $H_{d,g}$ to an extremal \c \ of the same degree
and genus.
\begin{enumerate}
\item \ Any disjoint union of lines \cite{H9}.
\item \ Any smooth \c \ with $d\geq g+3$ \cite{H9}.
\item \ Any ACM \c \ \cite{H9}.
\item \ Any \c \ in the biliaison equivalence class of an
extremal \c \ \cite{S}.
\item \ Any \c \ whose Rao module is a complete intersection
(also called a Koszul module). \cite{P}.
\end{enumerate}

\section{Techniques used}

The classification part of the problem uses standard methods.
What is new in studying connectedness questions is to prove the
existence of flat families of \c s, whose general member lies in one
\irr component, and whose special member lies in another \irr component.
We discuss here the different methods used to construct such families.

\subsection{Explicit equations}
If one knows the equations of the two types of \c s, one can attempt to
make a flat family by writing equations depending on a parameter $t$.
A simple example of this is the family of twisted cubic \c s having as
limit a plane nodal \c \ with an \emb ded point \cite[III, 9.8.4]{H2}.
Explicit equations are used in the papers
\cite{N2},\cite{N3},\cite{H10}. This technique is obviously
limited to situations where one has only to deal with very explicit
examples of \c s.

\subsection{Line drawings}
This method, used in \cite{H9}, is an extension of the first.
Some families of multiple structures on lines are proved by
explicit equations. These are then used as lemmas in drawings of much
more complicated \c s, supported on unions of lines. Combined with the
complete description of \c s contained in a double plane [15], this
allows one to show existence of families for many types of smooth \c s
that specialize to stick figures, such as the nonspecial \c s with
$d\geq g+3$.

\subsection{Triades}
This method is the most sophisticated, and potentially the most
powerful, but also the most technically difficult.
This method is developed in the three papers
\cite{HMP1},\cite{HMP2},\cite{HMP3} and applied in the papers
\cite{A1},\cite{A2},\cite{P}.

The idea is to develop an \alg theory of flat families, biliaison, and
Rao modules similar to the well known theory of biliaison and Rao module
for individual \c s \cite{MP2}.

In a flat family $C_t$ of \c s in ${\P}^3$ parametrized by a parameter
scheme $T$, the Rao module $M(C_t)$ is not in general constant in the
family. Also the sheaf analogue $\bigoplus_{n\in{\z}} R^1f_*({\mathcal
I}_C(n))$ as a sheaf of graded $S$-modules over $T$, does not commute
with base extension. So, for example, if $T={\mb{Spec}}\, A$ is affine,
one is led to consider the functor on $A$-modules
\[
M\mapsto \bigoplus_{n\in{\z}} H^1({\mathcal I}_C(n)\otimes_A M) \ .
\]
This is a coherent functor in the sense of Auslander \cite{H8},
but it is still not a fine enough invariant to play the role of the Rao
module for a family. So instead we consider the {\it triad} associated
to the family $C$: it is a 3-term complex $L_1\to L_0\to L_{-1}$ of
graded $S_A$-modules, where $S_A=A[x_0,x_1,x_2,x_3]$ is the homogeneous
coordinate ring of ${\P}^3_A$, whose middle cohomology retrieves the
cohomology  $H^1({\mathcal I}_C(n))$ of the family, and which satisfies
certain other technical conditions (see \cite[1.10]{HMP1} for the precise
definition). There is a notion of pseudoisomorphism for triades
\cite[1.7]{HMP1}, and then one obtains the analogue of Rao's theorem,
that two families of \c s are in the same biliaison equivalence class if
and only if their triades are  pseudoisomorphic up to shift in degrees
\cite[3.9]{HMP1}.

There is also an algorithmic method of constructing the universal family
of \c s associated to a triad \cite{HMP3}, and this becomes the basic
method of constructing flat families of \c s. The difficulty is that the
triad is not determined simply by knowing the Rao modules of the general
\c \ and the special curve: there are other choices to be made to
determine the triad. Thus to show the existence of a family connecting
\c s of particular types, one has to choose carefully a suitable triad
to give the family. This means also that while the method of triads is
good for making families, it is more difficult to prove the
non-existence of families between given types of \c s.

See also \cite{H6} for a slightly less brief introduction to the theory
of triades.

\section{An example}

Here we describe an example for which it is not yet known whether
$H_{d,g}$ is connected or not.

We consider smooth \c s $C$ of bidegree (3,7) on a smooth quadric
surface $Q$ in ${\P}^3$. Then $d=10$, $g=12$. We do not know if these \c
s can be connected to extremal \c s. Because of semicontinuity, these \c
s cannot be specializations of a family of \c s not contained in quadric
surfaces. So these \c s form an open subset of an \irr component of
$H_{10,12}$. The only possibility for connecting them to other \c s
requires specializing the quadric surface $Q$ to a quadric cone, the
union of two planes, or a double plane. One can show that if $Q$
specializes to a cone or to a union of two planes, the \c s must
necessarily acquire \emb ded points \cite[$\S$2]{H7}. So the only case
remaining is when $Q$ specializes to a double plane. Since one knows all
about \c s in the double plane \cite{H10}, it would be
sufficient to show the existence of a flat family going from
the \c s $C$ to a locally \CM
\c \ in the double plane, but this question has so far resisted analysis.

Another approach is to use biliaison. If one has a flat family going
from a \c \ $C_0$ to an extremal \c \ $E_0$, then by biliaison of the
family one obtains a flat family from $C_1$ to $E_1$, where $C_1$ and
$E_1$ are in the biliaison classes of $C_0$ and $E_0$, respectively.
Schlesinger's result \cite{S} shows that $E_1$ can be connected to an
extremal \c \ with the same degree and genus.

Now our curve $C$ of bidegree (3,7) on $Q$ is in the biliaison class of
a \c \ $C_0$ consisting of four skew lines, and one knows that four skew
lines can be connected to an extremal \c \ (cf.~$\S$2E above). The catch
is that in order to perform a biliaison of the family on the quadric
surface, the entire family must be contained in quadric surfaces.
In the case of four skew lines on $Q$, we do not know if they can be
specialized {\it on a quadric surface} to an extremal curve. The way we
know they are connected to an extremal \c \ is to pull them off the
quadric surface, giving more room to move around and then specialize.

So it seems that this example is a test case for the connectedness
question, and might possibly lead to a counterexample.

\setlength{\baselineskip}{16pt}


\begin{thebibliography}{10}
\addcontentsline{toc}{section}{References}
\bibitem[1]{A1} A\"it-Amrane, S. Sur le sch\'ema de Hilbert des
courbes gauches de degr\'e d et genre $g=(d-3)(d-4)/2$.
Th\`ese, Universit\'e Paris--Sud (1998).

\bibitem[2]{A2} A\"it-Amrane, S., Perrin, D. Un contre-exemple sur les
familles de courbes gauches. {\it Comm. Alg.} {\bf 28} (2000),
6003--6015.

\bibitem[3]{E} Ellingsrud, G. Sur le sch\'ema de Hilbert des
vari\'et\'es de codimension 2 dans ${\P}^e$ \`a c\^one de \CM.
{\it Ann. Sc. Ec. Norm. Sup} {\bf 8} (1975), 423--432.

  \bibitem[4]{GP} Gruson, L., Peskine, C. Genre des courbes de l'espace
projectif, in: {\it Algebraic Geometry}, (Troms{\o}, 1977),
Springer Lecture Notes in Math {\bf 687} (1978), 31--59.

\bibitem[5]{Ha} Halphen, G. M\'emoire sur la classification des
courbes gauches alg\'ebriques. {\it J. Ec. Polyt.} {\bf 52} (1882),
1--200.

\bibitem[6]{H1} Hartshorne, R. Connectedness of the Hilbert scheme.
{\it Publ. Math. IHES} {\bf 29} (1966), 5--48.


\bibitem[7]{H2} ---\!\!\!---\!\!\!------\!\!\!---\!\!\! \quad   {\it
Algebraic Geometry}, Springer (1977).

\bibitem[8]{H3} ---\!\!\!---\!\!\!------\!\!\!---\!\!\! \quad   Genre des
courbes alg\'ebriques dans l'espace projectif. {\it Sem. Bourbaki} {\bf
592} (1981/82).

\bibitem[9]{H4} ---\!\!\!---\!\!\!------\!\!\!---\!\!\!  \quad  The
genus of space curves. {\it Ann. Univ. Ferrara-Sez.VII - Sc. Mat.} {\bf
40} (1994), 207--223.

\bibitem[10]{H5}  ---\!\!\!---\!\!\!------\!\!\!---\!\!\!   \quad
Classification of algebraic space curves, III, in: {\it Algebraic
Geometry and Its Applications}, editor C.L.Bajaj (Springer, 1994),
113--120.

\bibitem[11]{H6} ---\!\!\!---\!\!\!------\!\!\!---\!\!\! \quad   Coherent
functors and families of space curves. {\it Rend. Sem. Mat. Fis. Milano}
{\bf 67} (1997), 87--93.


\bibitem[12]{H7} ---\!\!\!---\!\!\!------\!\!\!---\!\!\! \quad
Families of curves in ${\P}^3$ and Zeuthen's problem. {\it Memoirs Amer.
Math. Soc.} {\bf 130} (1997), no.~617.

\bibitem[13]{H8} ---\!\!\!---\!\!\!------\!\!\!---\!\!\! \quad   Coherent
functors.  {\it Advances in Math.} {\bf 140} (1998), 44--94.


\bibitem[14]{H9} ---\!\!\!---\!\!\!------\!\!\!---\!\!\! \quad On the
connectedness of the Hilbert scheme of \c s in ${\P}^3$. {\it Comm.
Alg.} {\bf 28} (2000), 6059--6077.

\bibitem[15]{H10} Hartshorne, R., Schlesinger, E. Curves in the double
plane. {\it Comm. Alg.} {\bf 28} (2000), 5655--5676.

\bibitem[16]{HMP1} Hartshorne, R., Martin--Deschamps, M., Perrin, D.
Triades et familles de courbes gauches. {\it Math. Ann.} {\bf 315}
(1999), 397--468.

\bibitem[17]{HMP3}
---\!\!\!---\!\!\!------\!\!\!---\!\!\!  \quad
Construction de familles de courbes gauches.
{\it Pacific J. Math.} {\bf 194} (2000), 97--116.

\bibitem[18]{HMP2}
---\!\!\!---\!\!\!------\!\!\!---\!\!\! \quad
Un th\'eor\`eme de Rao pour les familles de courbes gauches.
{\it J. Pure Appl. Algebra} {\bf 155} (2001), 53--76.


\bibitem[19]{Mac} Macaulay, F.S. {\it The Algebraic Theory of Modular
Systems}\, (Cambridge University Press, 1916).

\bibitem[20]{MP2}
  Martin--Deschamps, M.,
Perrin, D.  Sur la classification des courbes gauches. {\it
Ast\'erisque} {\bf 184--185} (1990).

\bibitem[21]{MP4}
---\!\!\!---\!\!\!------\!\!\!---\!\!\!  \quad  Sur les bornes du module
de Rao. {\it C.R.Acad. Sci. Paris} {\bf 317} (1993), 1159--1162.

\bibitem[22]{MP3}
---\!\!\!---\!\!\!------\!\!\!---\!\!\!  \quad  Le
sch\'ema de Hilbert des courbes localement de \CM n'est (presque) jamais
r\'eduit.
{\it Ann. Sc. Ec. Norm. Sup.} {\bf 29} (1996), 757--785.


\bibitem[23]{MP1} ---\!\!\!---\!\!\!------\!\!\!---\!\!\!  \quad Triades
et d\'eformations de sons-quotients. {\it Comm. Alg.} {\bf 28} (2000),
5601--5611.


\bibitem[24]{N1} Nollet, S. Subextremal curves. {\it manusc. math.}
{\bf 94} (1997), 303--317.

\bibitem[25]{N2} ---\!\!\!---\!\!\!------\!\!\!---\!\!\!  \quad
The Hilbert scheme of degree three curves. {\it Ann. Scient. ENS}
{\bf 30} (1997),\\ 367--384.

\bibitem[26]{N3} ---\!\!\!---\!\!\!------\!\!\!---\!\!\!  \quad
A remark on connectedness in Hilbert schemes.
{\it Comm. Alg.} {\bf 28} (2000),\\ 5745--5747.

\bibitem[27]{NS} Nollet, S., Schlesinger, E. Curves of degree four (to
appear).

\bibitem[28]{P} Perrin, D. Un pas vers la connexit\'e du sch\'ema de
Hilbert: les courbes de Koszul sont dans la composante des extr\'emales
(preprint).

\bibitem[29]{Sa} Sabadini, I. On the Hilbert scheme of curves of
degree $d$ and genus $\left(\begin{array}{c}d-3\\2\end{array}\right)-1$,
preprint.

\bibitem[30]{S} Schlesinger, E. Footnote to a paper by Hartshorne.
{\it Comm. Alg.} {\bf 28} (2000), 6079--6083.

\end{thebibliography}
\end{document}